\documentclass[11pt,reqno]{amsart} 
 \usepackage{amsmath,amssymb,amsthm}
 \usepackage{graphicx}
 \usepackage{cite}
 \usepackage[mathscr]{eucal}

 \theoremstyle{plain}
 \newtheorem{theorem}{Theorem}[section]  
 \newtheorem{corollary}[theorem]{Corollary}  
 \newtheorem{lemma}[theorem]{Lemma}  
 \newtheorem{proposition}[theorem]{Proposition}
 \newtheorem*{theorem*}{Theorem}

 \theoremstyle{plain}  
 
 \newtheorem{definition}[theorem]{Definition}

 \newtheoremstyle{citing}
   {3pt}
   {3pt}
   {\itshape}
   {}
   {\bfseries}
   {.}
   {.5em}
   {\thmnote{#3}}

 \theoremstyle{citing}

\numberwithin{equation}{section}
\allowdisplaybreaks[1]

\newlength{\intwidth}
\makeatletter
\DeclareRobustCommand{\cpvint}[2]
    {\mathop{%
       \text{%
         \settowidth{\intwidth}{%
           \ifx\ilimits@\displaylimits
             $\int_{#1}^{#2}$%
           \else
             $\int$%
           \fi}%
         \makebox[0pt][l]{\makebox[\intwidth]{$\text{C}$}}%
         $\int_{#1}^{#2}$}}}
\makeatother

\makeatletter
\DeclareRobustCommand{\cpvintsmall}[2]
    {\mathop{%
       \text{%
         \settowidth{\intwidth}{%
           \ifx\ilimits@\displaylimits
             $\int_{#1}^{#2}$%
           \else
             $\int$%
           \fi}%
         \makebox[0pt][l]{\makebox[\intwidth]{$\text{{\tiny C}}$}}%
         $\int_{#1}^{#2}$}}}
\makeatother

\newcommand{\rand}{\partial} 
\newcommand{\where}{:\:}

\newcommand{\laplace}{\Delta}

\newcommand{\nz}{{\mathbb N}}

\newcommand{\rz}{{\mathbb R}}  
\newcommand{\zz}{{\mathbb Z}}  
 
\newcommand{\eps}{\varepsilon}  
\renewcommand{\phi}{\varphi} 
\newcommand{\eval}{\vert}

\begin{document}
 
\title[Closed magnetic geodesics] {Alexandrov embedded closed magnetic geodesics on $S^2$}
\author{Matthias Schneider}
\address{Ruprecht-Karls-Universit\"at\\
         Im Neuenheimer Feld 288\\
         69120 Heidelberg, Germany\\}
\email{mschneid@mathi.uni-heidelberg.de} 
\date{April 8, 2010}  

\begin{abstract}
We prove the existence of two Alexandrov embedded closed magnetic geodesics on any
two dimensional sphere with nonnegative Gau{\ss} curvature.
\end{abstract}

\maketitle

\section{Introduction}
\label{sec:introduction}
Let $(S^2,g)$ be the two dimensional sphere equipped with a smooth metric $g$ and
$k:S^2 \to \rz$ a smooth positive function. We consider the following two equations
for curves $\gamma$ on $S^2$:
\begin{align}
\label{eq:magnetic_geodesic}
D_{t,g} \dot \gamma = k(\gamma) J_{g}(\gamma) \dot \gamma,
\end{align}
and
\begin{align}
\label{eq:1}
D_{t,g} \dot \gamma = |\dot \gamma|_{g} k(\gamma) J_{g}(\gamma) \dot \gamma,
\end{align}
where $D_{t,g}$ is the covariant derivative with respect to $g$,
and $J_g(x)$ is the rotation by $\pi/2$ in $T_xS^2$ measured with $g$.\\
Equation (\ref{eq:magnetic_geodesic}) describes the motion
of a charge in a magnetic field corresponding to the magnetic form
$k dV_g$
and solutions to \eqref{eq:magnetic_geodesic}
will be called {\em ($k$-)magnetic geodesics}
(see \cite{MR890489,MR1432462,MR2036336}).
Equation (\ref{eq:1}) corresponds to the problem of {\em prescribing geodesic curvature}, 
as its solutions $\gamma$ 
are constant speed curves with geodesic curvature $k_g(\gamma,t)$ given by $k(\gamma(t))$.\\
We study the existence of closed curves with prescribed geodesic curvature 
or equivalently the existence of periodic magnetic geodesics on prescribed energy levels 
$E_c:= \{(x,V) \in TS^2\where |V|_g=c\}$.\\
For fixed $k$ and $c>0$ 
the equations \eqref{eq:magnetic_geodesic} and \eqref{eq:1} are equivalent in the following
sense: If $\gamma$ is a nonconstant solution of \eqref{eq:1} with $k$ replaced by
$k/c$, then the curve 
$\gamma_c(t):=\gamma(ct/|\dot \gamma|_g)$
is a $k$-magnetic geodesic in $E_c$, and a $k$-magnetic geodesic in $E_c$
solves \eqref{eq:1} with $k$ replaced by $k/c$.\\
Solutions to \eqref{eq:magnetic_geodesic} or \eqref{eq:1} 
are invariant under a circle action:
For $\theta \in S^{1}=\rz/\zz$ and a closed curve $\gamma$
we define a closed curve $\theta*\gamma$ by
\begin{align*}
\theta*\gamma(t) = \gamma(t+\theta).  
\end{align*}
Thus, any solution gives rise to a $S^{1}$-orbit of solutions and we say
that two solutions $\gamma_1$ and $\gamma_2$ are (geometrically) distinct, if
$S^{1}*\gamma_1 \neq S^{1}*\gamma_2$.\\
There are different approaches to this problem, the Morse-Novikov theory
for (possibly multi-valued) variational functionals (see \cite{MR1185286,MR730159,MR1133303}), the theory of dynamical
systems using methods
from symplectic geometry (see \cite{MR902290,MR890489,MR1432462,MR1417851,MR2208800}) 
and Aubry-Mather's theory (see \cite{MR2036336}), and recently degree theory 
for immersed closed curves (see \cite{arXiv:0808.4038,2010arXiv1010.1879R}).\\
We follow \cite{arXiv:0808.4038} and consider solutions to \eqref{eq:1} as zeros of 
the vector field $X_{k,g}$ defined on the Sobolev space $H^{2,2}(S^1,S^2)$ as follows:
For $\gamma \in H^{2,2}(S^1,S^2)$ we let $X_{k,g}(\gamma)$ be the unique weak solution of   
\begin{align}
\label{eq:def_vector_field}
\big(-D_{t,g}^{2} + 1\big)X_{k,g}(\gamma)= 
-D_{t,g} \dot\gamma + |\dot \gamma|_{g} k(\gamma)J_{g}(\gamma)\dot\gamma   
\end{align}
in $T_\gamma H^{2,2}(S^1,S^2)$.  
The uniqueness implies that any zero of $X_{k,g}$ is a weak solution of \eqref{eq:1}
which is a classical solution in $C^{2}(S^1,S^2)$ applying standard regularity
theory.\\ 
For $(S^2,g)$ and $k$ positive 
it is conjectured (see \cite[1994-35,1996-18]{MR2078115}), 
\begin{align}
\label{eq:conj_arnold}
\text{ 
$\forall c>0$ the set $E_c$ contains (at least) two closed $k$-magnetic geodesics,
}  
\end{align}
which is is true for small energy levels depending
on $g$ and $k$ (see \cite{MR1432462,MR1417851}). 
In \cite{arXiv:0808.4038} it is shown that \eqref{eq:conj_arnold} is true, if the metric $g$ is 
$\frac14$-pinched, i.e. the Gau{\ss} curvature $K_g$ satisfies
\begin{align*}
\sup K_g < 4\inf K_g.   
\end{align*}
In fact, if $g$ is $\frac14$-pinched and $k$ is a positive function,
then every positive energy level $E_c$ contains two embedded (simple) 
closed $k$-magnetic geodesics.
We shall extend the above existence results.  
Instead of working in the class of embedded curves we consider solutions, that
are oriented Alexandrov embedded.
\begin{definition}
\label{def:alexandrov}
Let $B \subset \rz^2$ denote the open ball of radius $1$ centered at $0\in \rz^2$.
An immersion $\gamma\in C^1(\rand B,S^2)$ will be called {\em oriented Alexandrov embedded},
if there is an immersion $F\in C^1(\overline{B},S^2)$, such that
$F\eval_{\rand B}=\gamma$ and $F$ is orientation preserving in the sense that
\begin{align*}
\langle DF\eval_{x}x, J_g(\gamma(x))\dot\gamma(x)\rangle_{T_{\gamma(x)}S^2,g}>0
\end{align*}
for all $x \in \rand B$.
\end{definition}
Usually, Alexandrov embedded curves are defined to be the boundary of immersed manifolds.
We restrict ourselves to the case of the ball $B$ and oriented immersions $F$.  
If we equip $B$ with the metric $F^*g$ induced by $F$, then the outer normal $N_B(x)$
at $x \in \rand B$ with respect to $F^*g$ satisfies
\begin{align*}
DF\eval_x N_B(x) = N_\gamma(x)
\end{align*}
where $N_\gamma(x)$ denotes the normal to the curve $\gamma$ at $x \in \rand B$
defined by
\begin{align*}
N_\gamma(x):= |\dot\gamma(x)|^{-1}J_g(\gamma(x))\dot\gamma(x).  
\end{align*}
We shall prove
\begin{theorem}
\label{thm_existence}
Let $g$ be a smooth metric on $S^2$ with nonnegative Gau{\ss} curvature
and $k$ a positive smooth function.
Then there are (at least) two oriented Alexandrov embedded curves in $C^2(S^1,S^2)$ that
solve \eqref{eq:1}.
\end{theorem}
After completing the present work, we were informed
that the existence of one Alexandrov embedded closed curve solving
\eqref{eq:1} for a constant function $k$ is shown in \cite{robeday} (unpublished, 
see also \cite{2010arXiv1010.1879R}).\\
The equivalence between (\ref{eq:magnetic_geodesic}) and (\ref{eq:1}) leads to
\begin{corollary}
\label{cor_existence}
Let $g$ be a smooth metric on $S^2$ with nonnegative Gau{\ss} curvature
and $k$ a positive smooth function.
Then every energy level $E_c$ contains 
(at least) two oriented Alexandrov embedded closed magnetic geodesics.
\end{corollary}
The proof of our existence results is organized as follows.
In Section \ref{s:preliminiaries} we set up notation
and recall the basic properties of the $S^1$-equivariant Poincar\'{e}-Hopf index, 
$$\chi_{S^1}(X_{k,g},M_A) \in \zz,$$ 
where
$M_A$ is the set of oriented Alexandrov embedded regular curves
in $H^{2,2}(S^1,S^2)$. 
For positive constants $k_0$
we shall show that
$$\chi_{S^1}(X_{k_0,g_{can}},M_A)=-2.$$
Section \ref{sec:apriori-estimate} contains the apriori estimate which implies
that the set of solutions to \eqref{eq:1} is compact in $M_A$,
if the Gau{\ss} curvature of $(S^2,g)$ is nonnegative.
The homotopy invariance of the $S^{1}$-equivariant Poincar\'{e}-Hopf index then leads
to the identity
\begin{align*}
\chi_{S^1}(X_{k,g},M_A)=\chi_{S^1}(X_{k_0,g_{can}},M_A)=-2.  
\end{align*}
The resulting proof of Theorem \ref{thm_existence} is given in Section \ref{sec:existence}.

\section{Preliminaries}
\label{s:preliminiaries}
Let $S^{2}=\rand B_1(0) \subset \rz^{3}$ be the standard sphere and 
$g_{can}$ be the induced round metric. 
We consider the set of Sobolev functions
\begin{align*}
H^{2,2}(S^{1},S^{2}) := \{\gamma \in H^{2,2}(S^{1},\rz^{3}) \where \gamma(t) \in S^2 
\text{ for a.e. } t \in S^{1}\},   
\end{align*}
which is a sub-manifold of the Hilbert space $H^{2,2}(S^{1},\rz^{3})$
contained in $C^{1}(S^{1},\rz^{3})$.
The tangent space  
at $\gamma \in H^{2,2}(S^{1}, S^{2})$ is given by
\begin{align*}
T_\gamma H^{2,2}(S^{1}, S^{2}) := 
\{V \in H^{2,2}(S^{1},\rz^{3}\where V(t) \in T_{\gamma(t)}S^{2} \text{ for all }t \in S^{1}\}.  
\end{align*}
A metric $g$ on $S^{2}$ induces a metric on $H^{2,2}(S^{1},S^{2})$  
by 
\begin{align*}
\langle W,V\rangle_{T_\gamma H^{2,2}(S^{1}, S^{2}),g} 
:= \int_{S^{1}}\Big\langle &\big(-(D_{t,g})^{2}+1\big)V(t),\\
&\big(-(D_{t,g})^{2}+1\big)W(t)\Big\rangle_{\gamma(t),g}\, dt.   
\end{align*}
Since $g$ and $k$ are smooth, $X_{k,g}$ is a smooth vector field (see
\cite[Sec. 6]{arXiv:0808.4038,MR0464304}) 
on the set $H^{2,2}_{reg}(S^{1},S^{2})$ of regular curves,
\begin{align*}
H^{2,2}_{reg}(S^{1},S^{2}) := \{\gamma \in H^{2,2}(S^{1},S^{2})\where \dot\gamma(t) \neq 0 \text{ for all }t\in S^{1}\}. 
\end{align*}
Furthermore, we call $\gamma$ a prime curve, if
the isotropy group $$\{\theta \in S^1\where \theta*\gamma=\gamma\}$$ of $\gamma$ is trivial.\\
In \cite{arXiv:0808.4038} an integer valued $S^{1}$-degree, $\chi_{S^{1}}(X,M)$, is constructed for
a class of $S^{1}$-equivariant vector fields $X$ on open 
$S^1$-invariant subsets $M$ of prime curves in $H^{2,2}(S^{1},S^2)$. It is shown, that
\begin{align*}
\chi_{S^{1}}(X_{k,g},M) \in \zz  
\end{align*}
is defined, whenever $X_{k,g}$ is proper in $M$, i.e.
the set $\{\gamma \in M \where X(\gamma)=0 \}$ is compact,
and that the $S^1$-degree does not change under homotopies in the class
of proper vector fields.\\
We denote by $M_A\subset H^{2,2}(S^1,S^2)$ the set
\begin{align*}
M_A:= \{\gamma \in H^{2,2}_{reg}(S^1,S^2) 
\where \gamma \text{ is oriented Alexandrov embedded.}\}.
\end{align*}
By Lemma \ref{lem:alexandrov_open} below, $M_A$ is an open set of prime curves.\\
Equation \eqref{eq:1} with $g=g_{can}$ and $k \equiv k_0>0$ is analyzed in 
\cite[Sec. 5]{arXiv:0808.4038}:
The zeros of $X_{k_0,g_{can}}$ are given by $n$-fold iterates of a $S^2$-family of 
simple curves, and
there holds
\begin{align*}
\chi_{S^1}(X_{k_0,g_{can}},M_0)=-2,
\end{align*}
where $M_0$ denotes the set of embedded curves in $H^{2,2}_{reg}(S^1,S^2)$.
Since $n$-fold iterates are oriented Alexandrov embedded,
if and only if $n=1$ by Lemma \ref{lem:n_fold_alexandrov} below, the zeros of $X_{k_0,g_{can}}$
in $M_A$ and $M_0$ coincide and thus
\begin{align}
\label{eq:degree_unperturbed}
\chi_{S^1}(X_{k_0,g_{can}},M_A)=-2.  
\end{align}

\section{Alexandrov embedded curves}
\label{sec:stable}
We collect some properties of Alexandrov embedded curves known
to the experts, but difficult to find in the literature. 
For the sake of the readers convenience we include the proofs of these facts. 
\begin{lemma}
\label{lem:n_fold_alexandrov}
For $n \in \nz$ let $\gamma_n$ be the $n$-fold iterate of a
curve $\gamma_1 \in C^1(S^1,S^2)$ in $(S^2,g_{can})$ with nonnegative geodesic curvature.
If $\gamma_n$ is oriented Alexandrov embedded, then necessarily $n=1$. 
If, moreover, $\gamma_1$ is simple, 
then  $\gamma_n$ is oriented Alexandrov embedded if and only if $n=1$.     
\end{lemma}
\begin{proof}
If $\gamma_1$ is simple, $S^2\setminus \gamma_1(S^1)$ consists of two simply connected components
by the Jordan curve theorem.
We let $\tilde{B}$ be the component, where the normal $N_\gamma$ to $\gamma$
is the outer normal. By the Riemann mapping theorem there is an orientation preserving
diffeomorphism from $\tilde{B}$ to the open ball $B$ showing that
$\gamma_1$ is oriented Alexandrov embedded.\\
Assume $\gamma_n$ is oriented Alexandrov embedded and let $F_n$ be the 
corresponding immersion. We equip $B$ with the metric $F_n^*g_{can}$ induced by $F_n$. 
To obtain a contradiction we assume that $F_n$ is surjective. From
the Gau{\ss}-Bonnet formula theorem we derive
\begin{align*}
2\pi &= \int_{\rand B} k_{F_n*g_{can}}\, dS_{F_n*g_{can}} + \int_B K_{F_n^*g_{can}}\, dA_{F_n^*g_{can}}\\
&= \int_{\gamma_n} k_{g_{can}} \, ds_{g_{can}} + \int_B K_{F_n^*g_{can}} \, dA_{F_n^*g_{can}}\\
&\ge \int_{F_n(B)} K_{g_{can}}\, dA_{g_{can}} \ge \int_{S^2} K_{g_{can}} \, dA_{g_{can}}\\
&= 4\pi,
\end{align*}
which leads to the desired contradiction. Hence $F_n$ is not surjective and 
using stereographic coordinates
we may assume 
that $\gamma_1$ is a curve in the plane $(\rz^2,\delta)$ with standard metric $\delta$.
If we apply the Gau{\ss}-Bonnet formula to $(B,F_n^*\delta)$ and the
curve $\gamma_1$ in the plane, we obtain using the rotation index $i_{\gamma_1}$ 
\begin{align*}
2\pi &= \int_{\rand B} k_{F_n^*\delta}\, dS_{F_n^*\delta} + \int_B K_{F_n^*\delta}\, dA_{F_n^*\delta}\\ 
&= \int_{\gamma_n} k_{\delta} \, dS_{\delta}
= n\int_{\gamma_1} k_{\delta} \, dS_{\delta} = n i_{\gamma_1} 2\pi,
\end{align*}
which is only possible for $n=1$.
\end{proof}

\begin{lemma}
\label{lem:alexandrov_stable}
Let $(\gamma_n)$ in $C^2(\rand B,S^2)$ be a sequence of immersions,
which are oriented Alexandrov embedded, such that $(\gamma_n)$ converges 
to an immersion $\gamma_0$ in $C^2(\rand B, S^2)$ with strictly positive
geodesic curvature.
Then $\gamma_0$ is oriented Alexandrov embedded.    
\end{lemma}
\begin{proof}
We fix $(\gamma_n)$ and $\gamma_0$ that satisfy the above assumptions
and denote by $F_n:\:\overline{B} \to S^2$ the corresponding sequence of oriented immersions,
such that $F_n\eval_{\rand B}=\gamma_n$. As $\gamma_n$ is a $C^2$-map, we may
assume $F_n$ is in $C^2(\overline{B}, S^2)$ as well.\\
Since the convergence is in $C^2(\rand B, S^2)$, there is $\eps_0>0$
such that for all $n\in \nz_0$ the map
$T_{S^2,n}: \rand B\times (-\eps_0,\eps_0)\to S^2$ defined by
\begin{align*}
T_{S^2,n}(y,t) := Exp_{\gamma_n(y),g}\big(tN_{S^2,n}(y)\big)  
\end{align*}
is an immersion, where $N_{S^2,n}(y)$ is the outer normal at 
$\gamma_n(y)$ to $F_n(\rand B)$, if $n\neq 0$, and 
$N_{S^2,0}(y)=\lim_{n \to \infty}N_{S^2,n}(y)
=J_g(\gamma_0(y))|\dot\gamma_0(y)|^{-1} \dot\gamma_0(y)$. 
Moreover, we may assume that
the geodesic curvature 
of the curves $T_{S^2,n}(\cdot,s)$ is uniformly bounded, i.e. 
there is $k_0>0$ such that
\begin{align}
\label{eq:alex_tub_curv_bound}
k_0 \le  k_{g}\big(T_{S^2,n}(t,s),t\big) \le k_0^{-1}
\; \forall (n,t,s) \in \nz\times \rand B \times (-\eps_0,\eps_0). 
\end{align}
We consider $(\overline{B},F_n^*g)$, where $F_n^*g$ denotes the 
metric induced by $F_n$.
Shrinking $\eps_0>0$ we 
have the following 
\begin{proposition}
\label{prop:bound_inner_kragen}
For all $n\in \nz$ the map $T_{B,n}:\, \rand B \times (-\eps_0,0]
\to \overline{B}$,
\begin{align*}
T_{B,n}(y,t):= Exp_{y,F_n^*g}(t N_{B,n}(y)),  
\end{align*}
where $N_{B,n}(y)$ denotes the outer normal to $\rand B$ at $y$
with respect to $F_n^*g$, 
is well defined and a diffeomorphism onto its range.
\end{proposition}
\noindent
We postpone the proof of the proposition 
and proceed with the proof of Lemma \ref{lem:alexandrov_stable}.\\
Due to the unique solvability of the ordinary differential equation corresponding
to the exponential function we have
\begin{align}
\label{eq:2}
F_n\circ T_{B,n}(y,t)=T_{S^2,n}(y,t) \text{ for all }(y,t) \in \rand B \times (-\eps_0,0].   
\end{align}
Since $(\gamma_n)$ converges to $\gamma_0$ in $C^{2}(\rand B,S^2)$,
there is $n_0\in \nz$ such that a reparametrization of $\gamma_0$ 
is a graph over $\gamma_{n_0}$, i.e.
there is $s\in C^{1}(\rand B,(-\eps_0,\eps_0))$ 
and an oriented diffeomorphism $\alpha \in C^2(\rand B,\rand B)$, which 
is close to the identity,
such that
\begin{align*}
\gamma_0(\alpha(y)) = T_{S^2,n_0}(y,s(y))\quad \forall y \in \rand B.   
\end{align*}
If we define $F_0: \overline{B} \to S^2$ by
\begin{align*}
F_0(x) = 
\begin{cases}
F_{n_0}(x), &\text{if } x \text{ is in }\\
 &B \setminus T_{B,n_0}(\rand B \times (-\eps_0,0)),\\
T_{S^2,n_0}\big(y,(s(y)\eps_0^{-1}+1)t+s(y)\big),
&\text{if } x=T_{B,n_0}(y,t) \text{ for }\\
&(y,t)\in (\rand B \times (-\eps_0,0]), 
\end{cases}
\end{align*}
then $F_0 \in C^1(\overline{B},S^2)$, by \eqref{eq:2}. 
Since
$F_0\eval_{\rand B}=\gamma_0\circ \alpha$ and $F_0$ is 
an oriented immersion 
as a composition of such immersions, we see that $\gamma \circ \alpha$ is
oriented Alexandrov embedded. The diffeomorphism $\alpha$ can 
be extended to an oriented diffeomorphism $A$ of $\overline{B}$: We consider $\rand B$
as $\rz/2\pi \zz$, assume after a rotation $\alpha(0)=0$, and define $A$ in
polar coordinates by
\begin{align*}
A(r,\phi):= \big(r,\int_0^\phi s(r)+(1-s(r))\alpha'(\tau)\, d\tau\big),   
\end{align*}
where $s \in C^\infty([0,1],[0,1])$ is any function satisfying
\begin{align*}
s(r)=1 \text{, if } 0\le r \le \frac12 \text{ and }
s(r)=0 \text{, if }r=1  
\end{align*}
Consequently, $\gamma_0$ is oriented Alexandrov embedded as well using
$F_0 \circ A^{-1}$, which  
yields the claim of Lemma \ref{lem:alexandrov_stable}.\\
It remains to prove Proposition \ref{prop:bound_inner_kragen}:
Firstly, we note that for any $n \in \nz$ the map $T_{B,n}$ is defined 
in
\begin{align*}
U_n:= \{(y,t)\in \rand B \times (-\infty,0] \where -\sigma(y,n)\le t \le 0\},  
\end{align*}
where 
\begin{align*}
\sigma(y,n) &:= \sup\{t \where T_{B,n}(y,-t)\in \overline{B}\}
= \sup\{t \where T_{S^2,n}(y,-t)\in F_n(\overline{B})\},\\
\sigma(\rand B,n) &:= \inf\{\sigma(y,n)\where y \in \rand B\}.
\end{align*}
Differentiating \eqref{eq:2} we find in $U_n \cap \rand B \times (-\eps_0,0]$
\begin{align*}
DF_n\eval{T_{B,n}} \circ DT_{B,n}= DT_{S^2,n}.  
\end{align*}
Hence, $T_{B,n}$ is a local diffeomorphism and it is enough to show, 
after possibly shrinking $\eps_0>0$, that
\begin{align*}
\rand B \times (-\eps_0,0] \subset U_n \text{ for all }n \in \nz,\\
T_{B,n}\eval_{\rand B \times (-\eps_0,0]} \text{ is injective.}  
\end{align*}
To obtain a contradiction assume that
\begin{align*}
\sigma(\rand B,n)\to 0 \text{ as } n \to \infty.  
\end{align*}
It is standard to see that $\sigma(\rand B,n)$ is attained at some $y_{0,n}\in \rand B$
and that the geodesic 
\begin{align*}
[0,\sigma(\rand B,n)]\ni t \mapsto T_{B,n}(y,-t)  
\end{align*}
is perpendicular to $\rand B$ at $T_{B,n}(y,-\sigma(\rand B,n))$.\\ 
We fix $n_1>0$ such that $\sigma(\rand B,n_1)<\eps_0/2$ and 
$\sigma(\rand B,n_1)$ is attained at $y_0 \in \rand B$.
Due to the minimality of $\sigma(y_0,n_1)$ the parallel curve
\begin{align*}
y \mapsto T_{B,n_1}(y,\sigma(y_0,n_1))  
\end{align*}
lies inside $\overline{B}$ with positive curvature and touches
$\rand B$ at $y_0$ from the inside. This leads to the desired contradiction due to
the positive curvature of $\rand B$ and the maximum principle.\\
To show the injectivity we argue by contradiction
and assume
that there is a subsequence of $(\gamma_n)$, still denoted by $(\gamma_n)$,
and a sequence $(y_{1,n},y_{2,n},t_{1,n},t_{2,n})$
in $(\rand B)^2\times (0,\frac{1}n)^2$ such that
$y_{1,n}\neq y_{2,n}$ and
\begin{align*}
T_{B,n}(y_{1,n},-t_{1,n})=T_{B,n}(y_{2,n},-t_{2,n}).  
\end{align*}
Going to a subsequence we may assume
\begin{align*}
y_{1,n} \to y_1 \text{ and } y_{2,n} \to y_2
\text{ as } n \to \infty.    
\end{align*}
For $i\in \{1,2\}$ we have
\begin{align*}
F_n\circ T_{B,n}(y_{i,n},-t_{i,n})= T_{S^2,n}(y_{i,n},-t_{i,n}). 
\end{align*}
Consequently, $y_1\neq y_2$, because 
$(T_{S^2,n})$ converges in $C^1(\rand B\times (-\eps_0,\eps_0),S^2)$
to $T_{S^2,0}$, which is an immersion. 
By the same argument, we deduce that for any $0<r$ and any $0<\eps\le \eps_0$ 
there is $\delta=\delta(r,\eps)>0$ and $n_0=n_0(r,\eps)$
such that for all $n\ge n_0$
\begin{align*}
B_{\delta}(y_{1,n})\subset T_{S^2,n}(B_{r,\rand B}(y_{1,n}) \times (-\eps,\eps)).   
\end{align*}
Since $0<t_{1,n},t_{2,n}<\frac{1}{n}$ and $y_1\neq y_2$,
taking $0<r<\text{dist}_{\rand B}(y_1,y_2)/2$ and $0<\eps\le \eps_0$
we infer that $y_{2} \in B_{\delta}(y_{1})$ for $n$ large enough.
Consequently, as $\eps>0$ may be chosen arbitrarily small,
\begin{align*}
\sigma(\rand B,n)\to 0 \text{ as } n \to \infty.  
\end{align*}
This yields a contradiction as above and 
finishes the proof.
\end{proof}

\begin{lemma}
\label{lem:alexandrov_open}
The set of regular, oriented Alexandrov embedded curves is an open subset of prime curves
in $H^{2,2}(S^1,S^2)$.   
\end{lemma}
\begin{proof}
Let $\gamma_0 \in H^{2,2}(S^1,S^2)$ be oriented Alexandrov embedded. Then there
is an oriented immersion $F_0:\overline{B} \to S^2$ such that $F_0\eval_{\rand B}=\gamma$.
We may extend $F_0$ to an open neighborhood $U$ of $\overline{B}$ 
such that $F_0$ remains an immersion and $F_0(U)$ is open.
If $\gamma$ is close enough to $\gamma_0$, then $\gamma$ lies in $F_0(U)$ and 
we may write $\gamma$ as a graph
over $\gamma_0$ with respect to its normal. 
Since $\gamma(\rand B)\subset F_0(U)$
we may proceed as in the proof
of Lemma \ref{lem:alexandrov_stable} to deduce that $\gamma$ is
oriented Alexandrov embedded as well.\\
From Lemma \ref{lem:n_fold_alexandrov} we obtain that oriented Alexandrov embedded
curves are prime.   
\end{proof}

\section{The apriori estimate}
\label{sec:apriori-estimate}
We fix a continuous family of metrics $\{g_t\where t\in [0,1]\}$ on $S^2$
and a continuous family of positive continuous function $\{k_t \where t\in [0,1]\}$ on $S^2$,
such that the Gau{\ss} curvature $K_{g_t}$ is nonnegative and
\begin{align*}
k_{inf}:= \inf\{k_t(x)\where (x,t)\in S^2\times [0,1]\}>0.  
\end{align*}
We let $X_{t}$ be the vector field on $H^{2,2}(S^1,S^2)$ defined by
\begin{align*}
X_t:= X_{k_t,g_t}.  
\end{align*}
We shall show that the set 
\begin{align*}
X^{-1}(0) := \{(\gamma,t) \in M_A\times [0,1] \where X_t(\gamma)=0\}   
\end{align*}
is compact in $M_A\times [0,1]$.
Fix $(\gamma,t)\in X^{-1}(0)$. Then there is an oriented immersion $F:\overline{B}\to S^2$
with $F\eval_{\rand B}=\gamma$. We denote by $F^*g_t$ the induced metric on $B$.
Using $F\eval_{\rand B}=\gamma$, $K_{g_t}\ge 0$, and the fact that $F$ is a local isometry
from $(B,F^*g_t)$ to $(S^2,g_t)$, the Gau{\ss}-Bonnet formula gives
\begin{align*}
2\pi &=
\int_{\rand B}k_{F*g_t} \, dS_{F*g_t} +\int_{B}K_{F*g_t} \,dA_{F*g_t}\\
&\ge
\int_{\gamma}k_t \, dS_{g_t} \ge k_{inf} L(\gamma,g_t),
\end{align*}
where $L(\gamma,g_t)$ denotes the length of $\gamma$ in the metric $g_t$.\\
To obtain a contradiction assume that there is $(\gamma_n,t_n)$ in $X^{-1}(0)$
such that $L(\gamma_n,g_{t_n})\to 0$ as $n \to \infty$. We denote by $F_n$ the corresponding
oriented immersions. Since the sets $\{g_{t}:\: t\in [0,1]\}$ and
$\{k_{t}:\: t\in [0,1]\}$ are compact, 
all involved metric are uniform equivalent to the standard metric $g_{can}$ and
there is $C_{k,g_{can}}>0$, such that
\begin{align*}
|k_{g_{can}}(\gamma_n,t)| \le C_{k,g_{can}} \text{ for all } (n,t) \in \nz \times S^1.   
\end{align*}
For $n \in \nz$ we let
\begin{align*}
L_n &:= L(\gamma_n,g_{can}) = \int_{\rand B} dS_{F_n^*g_{can}},\\
\tilde{L}_n &:= \int_{\gamma_n}k_{g_{can}} dS_{g_{can}} 
= \int_{\rand B} k_{F_n^*g_{can}} dS_{F_n^*g_{can}},\\
A_n &:= \int_B dA_{F_n^*g_{can}} = \int_B K_{F_n^*g_{can}} dA_{F_n^*g_{can}}.   
\end{align*}
Due to the uniform equivalence of the involved metrics and the uniform bound on
the geodesic curvature
$L_n$ and $\tilde{L}_n$ tend to $0$ as $n \to \infty$. 
From the Gau{\ss}-Bonnet formula applied to $(B,F_n^*g_{can})$ we obtain
\begin{align*}
A_n = 2\pi - \tilde{L}_n.  
\end{align*}
Applying the isoperimetric inequality \cite{MR0500557} to $(B,F_n^*g_{can})$ we find
\begin{align*}
L_n^2 &\ge 4\pi A_n - (A_n)^2\\
&= 4\pi\big(2\pi - \tilde{L}_n\big) -\big(2\pi - \tilde{L}_n\big)^2\\
&= 4\pi^2- (\tilde{L}_n)^2, 
\end{align*}
which is impossible for large $n$.\\
Consequently, the length $L(\gamma,g_t)$ for $(\gamma,t)\in X^{-1}(0)$
satisfies
\begin{align}
\label{eq:length_bound}
c\le L(\gamma,g_t) \le \big(\inf\{k_t(x)\}\big)^{-1} 2\pi
\end{align}
for some positive constant $c=c(\{k_t\},\{g_t\})$.\\
Fix a sequence $(\gamma_n,t_n)_{n\in \nz}$ in $X^{-1}(0)$. Since $\gamma_n$ is a
zero of $X_{t_n}$, the curve $\gamma_n$ is parameterized proportional
to its arc-length. From the bound in \eqref{eq:length_bound} we see that
$(\gamma_n)$ is uniformly bounded in $C^1(S^1,S^2)$. Using the equation
\eqref{eq:1} and standard regularity theory we find that $(\gamma_n)$
is bounded in $C^3(S^1,S^2)$, such that we may extract a subsequence, still denoted
by $(\gamma_n,t_n)_{n \in \nz}$,
which
converges in $C^2(S^1,S^2)\times [0,1]$ to $(\gamma_0,t_0)$.
Due to the convergence in $C^2(S^1,S^2)$ we have $X_{t_0}(\gamma_0)=0$,
and the lower bound in \eqref{eq:length_bound} implies that $\gamma_0$ is an immersion.
From the stability of oriented Alexandrov embeddings in Lemma \ref{lem:alexandrov_stable}
we deduce that $\gamma_0$ is oriented Alexandrov embedded and hence
$(\gamma_0,t_0)\in X^{-1}(0)$. This shows that $X^{-1}(0)$ is compact.\\
From the homotopy invariance we now get 
\begin{align}
\label{eq:degree_equal_apriori}
\chi_{S^1}(X_{k_0,g_{can}},M_A)= \chi_{S^1}(X_{k_1,g_1},M_A).
\end{align}

\section{Existence results}
\label{sec:existence}
We give the proof of our main existence result.
\begin{proof}[{Proof of Theorem \ref{thm_existence}}]
From the uniformization theorem there are a function $\phi \in C^\infty(S^2,\rz)$
and an isometry $F$ of $(S^2,g)$
to $(S^2,e^\phi g_{can})$, where $g_{can}$ denotes the 
standard round metric on $S^2$. Since the problem of prescribing
geodesic curvature is invariant under isometries we may assume
without loss of generality that
\begin{align*}
g= e^\phi g_{can}. 
\end{align*}
We consider the family of metrics $\{g_t\where t\in [0,1]\}$ 
and the family of positive continuous function $\{k_t \where t\in [0,1]\}$
defined by
\begin{align*}
g_t &:= e^{t\phi} g_{can},\\
k_t &:= (1-t)(\inf k) +tk.  
\end{align*}
Then $k_t\ge \inf k>0$ for all $t\in [0,1]$ and
the Gau{\ss} curvature $K_{g_t}$ of the metric $g_t$ satisfies
\begin{align*}
K_{g_t} &= e^{-t\phi}\big(-t\laplace_{g_{can}}(\phi)+2\big)\\
&= e^{-t\phi}\big(-t(2-K_{g}e^\phi)+2\big)\ge 0,
\end{align*}
because $K_g$ is nonnegative.\\
From Section \ref{sec:apriori-estimate} the homotopy 
\begin{align*}
[0,1]\ni t\mapsto X_{k_t,g_t}  
\end{align*}
is $(M_A,g_t,S^1)$-admissible and by \eqref{eq:degree_unperturbed} 
and \eqref{eq:degree_equal_apriori} 
\begin{align*}
-2= \chi_{S^1}(X_{k_0,g_{can}},M_A)=\chi_{S^1}(X_{k,g},M_A).  
\end{align*}
Since the local degree of an isolated critical orbit
is larger than $-1$ by \cite[Lem 4.1]{arXiv:0808.4038}, there
are at least two solutions to \eqref{eq:1}.
This gives the claim.
\end{proof}

\section*{Acknowledgements}
I would like to thank Friedrich Tomi and Karsten Grosse-Brauckmann 
for valuable discussions and suggestions.

\bibliographystyle{plain}
\bibliography{geodesic_curves}

\end{document}